\documentclass[10pt]{amsart}
\usepackage{latexsym}
\usepackage{amsmath}
\usepackage{amssymb}
\usepackage[all]{xy}

\newtheorem{thm}{Theorem}[section]
\newtheorem{prop}[thm]{Proposition}
\newtheorem{que}[thm]{Question}
\newtheorem{cor}[thm]{Corollary}
\newtheorem{lem}[thm]{Lemma}
\theoremstyle{definition}
\newtheorem{defn}[thm]{Definition}

\newtheorem{rem}[thm]{Remark}
\newtheorem{ex}[thm]{Example}

\newcommand{\Q}{{\mathbb Q}}
\newcommand{\Z}{{\mathbb Z}}
\newcommand{\C}{{\mathbb C}}

\begin{document}

\title{  Pre-c-symplectic 
condition for  the  product  of odd-spheres  }
\author{Junro SATO \ and \ Toshihiro YAMAGUCHI}
\renewcommand{\thefootnote}{\fnsymbol{footnote}}
\footnote[0]{Our definition of {\it pre-c-symplectic} is completely 
different from usual one of 
{\it presymplectic}(cf. \cite{KT}, \cite{HW}).\\
MSC: 55P62,53D05\\
Keywords: symplectic, c-symplectic, pre-c-symplectic, Sullivan  model, rational homotopy type, almost free toral action, 
rational toral rank,  Hasse diagram of  rational toral ranks}
\date{}
\address{Faculty of Education, Kochi University, 2-5-1,Kochi,780-8520, JAPAN}
\email{junro@kochi-u.ac.jp, tyamag@kochi-u.ac.jp}
\maketitle

\begin{abstract}
We say that a simply connected space $X$ is  {\it pre-c-symplectic} 
if it  is the fibre of a rational fibration $X\to Y\to \C P^{\infty}$ where $Y$ is  cohomologically 
symplectic in the sense that
there is a degree 2 cohomology class which cups to a top class.
It is a rational homotopical property but
not a cohomological one.
By using Sullivan's minimal models \cite{FHT},
we give the necessary and sufficient condition that the  product  of odd-spheres 
$X=S^{k_1}\times \cdots \times S^{k_n}$ 
is pre-c-symplectic
and see some related topics.
Also  we give a charactarization of the Hasse diagram of   rational toral ranks
for  a space $X$ \cite{Y} 
as a
necessary condition  to be pre-c-symplectic
and see some examples in the cases of finite-oddly generated 
rational homotopy groups.
\end{abstract}

\section{Introduction}
Recall the question:{\it ``If a symplectic manifold is a nilpotent space, 
what special homotopical properties are apparent ? 
Conversely, what nilpotent spaces have symplectic or c-symplectic structure ?''}
\cite[4.99]{FOT}.  
Here a rationally Poincar\'{e} dual space  $Y$ 
(the 
graded algebra $H^*(Y;\Q)$ is a  Poincar\'{e} duality algebra\cite[Def.3.1]{FOT}) 
with formal dimension $$fd(Y):=\max \{i|H^i(Y;\Q )\neq  0\}$$
$=2n$ 
 is said to be {\it c-symplectic (cohomologically symplectic)}
if there is  a  rational cohomology 
class $\omega\in H^2(Y;\Q)$ such that $\omega^n$ is a top 
class for $Y$ \cite[Def.4.87]{FOT}(\cite{TO},\cite{McS}),
many of which are known to be realized by $2n$-dimensional 
smooth manifolds (\cite{FOT}).
A lot  of  results on the above problem and related topics 
are given in   rational homotopy theory
 (cf.\cite{LO1}, \cite{LO}, \cite{TO}, \cite{Ke}, \cite{FOT}, \cite{KM}, \cite{LM},  \cite{K}, \cite{BM}, 
\cite{BFM},..).  
For example,
G.Lupton and J.Oprea \cite{LO1}
study   the formalising tendency of
certain  symplectic manifolds
using techniques of D.Sullivan's rational model \cite{Su}.
Notice that
it is known  that  the connected sum $\C P^2\sharp \C P^2$ is c-symplectic but not symplectic
\cite{A}(\cite[page 263]{LO}), 
for the $n$-dimensional complex projective space $\C P^n$.
In \cite{Th}(\cite[Theorem 6.3]{McS}), \cite{Ke},  \cite{K}..., we can see 
conditions that a total space with a degree 2 cohomology class admits  
  a symplectic structure  
in  a certain  fibration.  
But we don't mention  anything about 
symplectic geometry in this paper.

For a simply connected c-symplectic space $Y$, we have
 $\omega \in {Hom(\pi_2(Y),\Q )}$ 
for the  class $\omega$ of $H^2(Y;\Q )$
from Hurewicz isomorphism.
In particular, $\pi_2(Y)\otimes \Q\neq 0$. 
So
there is a simply connected space $X$ 
 that is the fibre of a fibration $$X\to Y\to \C P^{\infty} \ \ \ (1)$$
where $\C P^{\infty}=\cup_{n=1}^{\infty} \C P^n(=K(\Z ,2))$, $\pi_*(X)\otimes \Q\oplus \Q\cdot  t^*=\pi_*(Y)\otimes \Q$ for a cohomology element $t$ with $\deg (t)=2$ (necessarily we don't need $t=\omega$) 
and $H^*(\C P^{\infty};\Q )=\Q [t]$.



\begin{defn}
 We say a simply connected space $X$ to be   {\it  pre-c-symplectic} 
({\it pre-cohomologically symplectic}) 
if $X$ is the fibre of a fibration (1) where $Y$ is  c-symplectic.
\end{defn}

For example,  $\C P^n$
is a symplectic manifold,  whose  pre-c-symplectic space must be
 the $2n+1$-dimensional sphere $S^{2n+1}$. 
 It is induced by the Hopf fibration
 $S^1\to S^{2n+1}\to \C P^n$ \cite[p.95]{Ar}.
We know that   $fd(Y)=2n$ if and only if  $fd(X)=2n+1$ 
in $(1)$
 from the Gysin exact sequence of of the induced fibration $S^1\to X\to Y$.
When $\dim \pi_2(Y)\otimes \Q >1$,
 $(1)$  may  not be rational homotopically unique 
for $Y$.
For example,
when $Y$ is $S^2\times  \C P^2$,
two spaces $S^3\times  \C P^2$
and $S^2\times  S^5$
are both 
its pre-c-symplectic spaces (there are three pre-c-symplectic 
spaces in  the case of \cite[Example 2.12]{LO1}).
The being c-symplectic   and the being pre-c-symplectic   are complementary.
If a  space is c-symplectic, it is not pre-c-symplectic
and moreover if a  space is pre-c-symplectic, it is not c-symplectic.
The being c-symplectic
is  preserved by  product; i.e.,
$Y_1\times Y_2$ is   pre-c-symplectic by the class $\omega_1+\omega_2$ 
when  $Y_1$ and $Y_2$ are both  c-symplectic by classes $\omega_1$ and 
$\omega_2$, respectively.
But the being pre-c-symplectic can not
since then the formal dimension is even.


Of course, the being  
  pre-c-symplectic depends on the rational  homotopy type of $X$.
Recall  the Sullivan's rational  model  theory \cite{Su}.
Let 
the  Sullivan minimal model of $X$ be $M(X)=(\Lambda {V},d)$. 
  It is a free $\Q$-commutative differential graded algebra (dga)
 with a $\Q$-graded vector space $V=\bigoplus_{i\geq 2}V^i$
 where $\dim V^i<\infty$ and a decomposable differential; i.e., $d(V^i) \subset (\Lambda^+{V} \cdot \Lambda^+{V})^{i+1}$ and $d \circ d=0$.
 Here  $\Lambda^+{V}$ is 
 the ideal of $\Lambda{V}$ generated by elements of positive degree. 
Denote the degree of a homogeneous element $f$ of a graded algebra as $|{f}|$.
Then  $xy=(-1)^{|{x}||{y}|}yx$ and $d(xy)=d(x)y+(-1)^{|{x}|}xd(y)$. 
Note that  $M(X)$ determines the rational homotopy type of $X$.
In particular,  it is known that  $$H^*(\Lambda {V},d)\cong H^*(X;\Q )
\mbox{ \ and \ } V^i\cong Hom(\pi_i(X),\Q).$$
Refer  \cite[\S 12$\sim$\S 15]{FHT} for detail.
Especially, 
$(1)$  is replaced   with the relative model (KS-model)  \cite{FHT}  
$$(\Q[t],0)\to (\Q[t]\otimes \Lambda V,D)\to (\Lambda V,d) \ \ \ \ \ (2)$$
where $|t|=2$
and $\overline{D}=d$.
We often say 
that $M(Y)=(\Q[t]\otimes \Lambda V,D)$
is c-symplectic when $Y$  is so.
When $\pi_*(X)\otimes \Q<\infty$ and $\dim H^*(X;\Q )<\infty$,
a simply connected space $X$ is said to be  elliptic.
It is known that
$$fd(X)=fd(\Lambda V,d)=
\sum_i|y_i|-\sum_i(|x_i|-1)$$
for $V^{odd}=\Q (y_i)_i$ and  $V^{even}=\Q (x_i)_i$ 
when $X$ is elliptic \cite[\S 32]{FHT}.
When is a simply connected 
space $X$     pre-c-symplectic ?
Notice that if a pure model  $M(Y)=(\Lambda {U},d_Y)$,
which satisfies $d_YU^{even}=0$ and $d_YU^{odd}\subset \Lambda U^{even}$,
is c-symplectic, then $\dim U^{even}=\dim U^{odd}$ \cite{LO1}.
For example,
any simply connected symplectic homogeneous 
space is a maximal rank homogeneous 
space
 \cite[Corollary 2.5]{LO1}.
So, from $(2)$,  it may be  natural to expect that    $\dim V^{even}=\dim V^{odd}-1$
if a pure model $M(X)=(\Lambda {V},d)$  is pre-c-symplectic.
But it is  false (cf. Theorem \ref{A}  below).
If anything, 
{``it is relatively easy to construct c-symplectic Sullivan minimal models''} 
(cf.\cite[Example 2.9]{LO1},
\cite[p.263]{LO}) 
and furthermore {\it pre-c-symplectic spaces exist everywhere.}
The latter  is  nearly true if we can suitably
change the  ratio of degrees of basis elements
of $V$ for $M(X)=(\Lambda {V},d)$.   
For example, for any even dimensional simply connected compact manifold $B$, 
the product space $X=B\times S^N$ for
 the $N$-dimensional sphere $S^N$ 
 is pre-c-symplectic for any odd integer $N$ with 
$N>\dim B$.
Indeed, we can put the model of $(2)$ as
$M(Y)=(\Q [t]\otimes \Lambda V\otimes \Lambda v,D)$ by  $$D(v)=\alpha \cdot t^{(N+1-\dim B)/2}-t^{(N+1)/2}
\mbox{ \ and\ \ \ } D(b)=d_B(b)$$
for $b\in M(B)=(\Lambda V,d_B)$,
 the fundamental class $[\alpha]$ of $H^*(B;\Q)$
and $M(S^N)=(\Lambda v,0)$ with $|v|=N$.
Then 
$$H^*(Y;\Q )=H^*(B;\Q )[t]/(\alpha \cdot t^{(N+1-\dim B)/2}-t^{(N+1)/2})$$ and 
$[t]^{(\dim B+N-1)/2}=[\alpha\cdot t^{(N-1)/2}]\neq 0$.
Since $fd (Y)=\dim B+N-1$,
we see $Y$ is c-symplectic, that is,  $X$ is pre-c-symplectic.
In general,  it seems difficult to find
the smallest $N$ such that $X$ is pre-c-symplectic.
This is a symbolic example in this paper.


We will study the conditions  of  spaces to be     pre-c-symplectic,
especially in the most rational homotopically simple case, that is, 
we suppose that a finite simply connected 
complex $X$ has the rational cohomology structure of
the exterior algebra over $\Q$:
$$H^*(X;\Q)\cong \Lambda (v_1,v_2,\cdots ,v_n)$$
with $1<|v_1|=k_1\leq |v_2|=k_2\leq \cdots \leq |v_n|=k_n$ all odd.
Then $X$ has the rational homotopy type of the n-product of simply connected odd-spheres:
$$X\simeq_{\Q} S^{k_1}\times S^{k_2}\times \cdots \times S^{k_n}\ \ \ \  \ k_i;\ \mbox{odd}$$
($\simeq_{\Q}$ means ``is rational
homotopy equivalent to'')  and the Sullivan minimal model is given by
$$M(X)\cong (\Lambda (v_1,v_2,\cdots ,v_n),0).$$
For example,  simply  connected   compact Lie groups of rank $n$ satisfy the condition
(H.Hopf).
In this case,   (2) is written as
$$(\Q[t],0)\to (\Q[t]\otimes \Lambda (v_1,v_2,\cdots ,v_n),D)\to (\Lambda (v_1,v_2,\cdots ,v_n),0).$$




In this paper, we show 

\begin{thm} \label{A}When 
$H^*(X;\Q)\cong \Lambda (v_1,v_2,\cdots ,v_n)$ with all $|v_i|$ odd
and  $1<|v_1|\leq |v_2|\leq \cdots \leq |v_n|$,
then $X$ is pre-c-symplectic 
if and only if $n$ is odd
and $|v_1|+|v_{n-1}|<|v_n|$,
 $|v_2|+|v_{n-2}|<|v_n|$, $\cdots$,
 $|v_{(n-1)/2}|+|v_{(n+1)/2}|<|v_n|$.
\end{thm}

\begin{rem}\label{r13} The ``{\it if}'' part   of 
Theorem \ref{A} does not follow 
when 
$H^*(X;\Q )$ is not free; i.e.,
$d\neq 0$ for
$M(X)=(\Lambda (v_1,\cdots ,v_n),d)$. 
For example, when $M(X)=(\Lambda (v_1,v_2,v_3,v_4,v_5),d)$ 
with  $|v_1|=3$, $|v_2|=|v_3|=5$, $|v_4|=9$, 
$|v_5|=13$, $dv_1=dv_2=dv_3=dv_5=0$
and  
$dv_4=v_2v_3$,
any 
model $(\Q[t]\otimes \Lambda (v_1,v_2,v_3,v_4,v_5),D)$ of $(2)$
is not  pre-c-symplectic.
Indeed,  
the element $v_1v_4$
can not be a $D$-cocycle and 
$Dv_5$ can not contain the cocycle $v_iv_4t$ for $i=2,3$ from degree reasons.
So we can not construct the form $Dv_5=v_av_bt^*+v_cv_dt^*+t^7$
with $\{ a,b,c,d\} =\{ 1,2,3,4\}$. 
Also the ``{\it only if}'' part   of 
Theorem \ref{A} does not follow 
when 
$H^*(X;\Q )$ is not free.
For example, when $n=3$, $|v_1|=|v_2|=3$, $|v_3|=5$, $dv_1=dv_2=0$, 
$dv_3=v_1v_2$,
 the model $(\Q[t]\otimes \Lambda (v_1,v_2,v_3),D)$ of $(2)$
 with $Dv_1=Dv_2=0$ and 
$Dv_3=v_1v_2+t^3$   
  is c-symplectic by $[t^5]\neq 0$ 
but $|v_1|+|v_2|>|v_3|$ (see Theorem \ref{odd}). 
 \end{rem}

\begin{cor}\label{L}
Let $X$ be a compact connected simple Lie group $G$  of rank$\ G>1$.
Then  $X$ is  pre-c-symplectic if and only if  
$G$ is  $B_n$ or $C_n$ with $n$ odd, 
or  $E_7$.
\end{cor}

For example,   for  the 5-th symplectic group $Sp(5)$,
the rational cohomology is given as 
$H^*(Sp(5);\Q)= \Lambda (v_1,v_2,v_3,v_4, v_5)$
with the degrees $|v_1|=3$,  $|v_2|=7$, $|v_3|=11$, $|v_4|=15$ and  $|v_5|=19$.
From Corollary \ref{L},  it is pre-c-symplectic. 
There are at least the  four rational homotopy types of  c-symplectic models: 
$$\begin{array}{ll}i)&\ \ Dv_5=v_1v_4t+v_2v_3t+t^{10}, \ Dv_1=Dv_2=Dv_3=Dv_4=0\\
ii)&\ \ Dv_5=v_1v_4t+v_2v_3t+t^{10}, \ Dv_3=v_1v_2t,\ Dv_4=0\\
iii)&\ \ Dv_5=v_1v_4t+v_2v_3t+t^{10}, \ Dv_3=0,\ Dv_4=v_1v_3t\\
iv)&\ \ Dv_5=v_1v_4t+v_2v_3t+t^{10}, \ Dv_3=v_1v_2t,\ Dv_4=v_1v_3t.
\end{array}$$
Although the cohomology algebra structures of them are very different,
 they are  all c-symplectic with formal dimension
$54$.
For example, the cohomology algebras of 
$i),\ ii)$ and $iv)$
are given  as\\
$i)\ \ \Q [t]\otimes \Lambda (v_1,v_2,v_3,v_4)/(v_1v_4t+v_2v_3t+t^{10})$\\
$ii)\ \ \Q [t,u_1,u_2]\otimes \Lambda (v_1,v_2,v_4)/(v_1v_4t+u_2t+t^{10},v_2u_1+v_1u_2,
v_1v_2t,v_1u_1,  v_2u_2)$\\
$iv)\ \ \Q [t,u_1,u_2,u_3]\otimes \Lambda (v_1,v_2)/$\\
\ \ \ \ \ \ \ \ $(u_2t+u_3t+t^{10},v_2u_1+v_1u_2,v_1v_2t, v_1u_1,  v_2u_2,v_1u_3,u_1u_2,u_1u_3,u_1t)$,\\
where $u_1=[v_1v_3]$, $u_2=[v_2v_3]$
and  $u_3=[v_1v_4]$.\\ 
\hspace{4mm}





Let $r_0(X)$ be  the {\it rational toral rank} of $X$,  which is 
 the largest integer $r$ such that an $r$-torus
 $T^r=S^1 \times\dots\times S^1$($r$-factors)  can act continuously
 on  a space $X'$ in  the rational homotopy type of $X$
 with all its isotropy subgroups finite (almost free action) \cite{H}, \cite{FOT}.
For example,
$r_0(S^{k_1}\times \cdots \times S^{k_n})=n$ when  $k_i$ are all odd
and $r_0(\C P^n)=0$.
Pre-c-symplectic spaces 
are  related to  almost free toral actions.
Indeed, for (1), 
there is
 a free $S^1$-action on a finite complex $X'$  with $X'_{\Q}\simeq X_{\Q}$,  
from S.Halperin's Proposition \ref{H} of \S 3.
Here $X_{\Q}$ means  the rationalization of $X$ \cite{HMR}.
Thus we have the Borel fibration
$$X'\to ES^1\times_{S^1}X'\to BS^1\ \ \ \ \ (3)$$
with $\dim H^*(ES^1\times_{S^1}X';\Q)<\infty$.
It is rationally equivalent to $(1)$. Namely,
\begin{thm}\label{Pre}
 A simply connected space $X$ is   pre-c-symplectic 
if and only if there is  rationally an almost free  circle action on 
 $X$ 
such that the orbit space
is c-symplectic.
\end{thm}

In particular,  we see that $r_0(X)>0$
for a pre-c-symplectic space $X$.
The being  c-symplectic is surely a cohomological property.
But the being pre-c-symplec depends 
on the dga  and not simply on its   cohomology.
For example,
when two spaces 
$X_1$ and $X_2$ are given by $X_1=(S^3\times S^8)\sharp (S^3\times S^8)$ and $M(X_2)=(\Lambda (v_1,v_2,v_3),d)$
with $|v_1|=|v_2|=3$, $|v_3|=5$, $dv_1=dv_2=0$ and  $dv_3=v_1v_2$, 
we have a graded algebra isomorphism
$$H^*(X_i;\Q)\cong \Lambda (x,y)\otimes \Q [w,u]/(xy, xu,xw+yu,yw,w^2,wu,u^2)$$
with $|x|=|y|=3$ and $|w|=|u|=8$
 for $i=1,2$. 
When $i=2$,  $u=[v_1v_3]$ and   $w=[v_2v_3]$.
 Recall that $r_0(X_1)=0$ \cite[Theorem 1.1(2)]{KY},
so $X_1$ can  not be pre-c-symplectic
from Theorem \ref{Pre},  but 
  $X_2$ is pre-c-symplectic (see Remark \ref{r13}).
The following proposition seems   a special case of 
\cite[Corollary 3.7, (Theorem 5.2)]{LO}.

\begin{prop}\label{SS}
For a  simply connected c-symplectic space $Y$, 
$r_0(Y)=0$.
\end{prop}

If   $ET^a\times_{T^a}^{\mu}X$ is c-symplectic
for some $T^a$-action $\mu$,
then ($ET^{a-1}\times_{T^{a-1}}^{\tau}X$  is pre-c-symplectic
for any restriction $\tau$ on $T^{a-1}$ of $\mu$ 
and)  
$ET^b\times_{T^b}^{\tau}X$  ($a\neq b$)  can   not be 
c-symplectic for  any  restriction or
 extension $\tau$ 
on $T^b$ of $\mu$
from Proposition \ref{SS}.
But notice that when   $X$ or  $ET^a\times_{T^a}^{\mu}X$ is pre-c-symplectic, 
 $ET^b\times_{T^b}^{\tau}X$  ($a<b$) may be  pre-c-symplectic for  an  
 extension $\tau$. 
It may   complicate the being  pre-c-symplectic than 
the being  c-symplectic.
For example,
when $X\simeq_{\Q}S^3\times S^3\times S^7$
with  $M(X)=( \Lambda (v_1,v_2 ,v_3),0)$,
$X$  
is  pre-c-symplectic
since the model 
$(\Q [t]\otimes \Lambda (v_1,v_2 ,v_3),D)$ of $(3)$
is given by  $Dv_1=Dv_2=0$ and  
$Dv_3=v_1v_2t+t^4$.
Indeed, then 
$fd(ES^1\times_{S^1}X)=12$ 
and $[t^6]\neq 0$ (See Example \ref{3}).
On the other hand,
for any almost free $T^2$-action on $X$,
the Borel space $ET^2\times_{T^2}X$ is  also  pre-c-symplectic
since  the model  of $(3)$
is given by Proposition \ref{H} as
$$(\Q [t_3],0)\to (\Q [t_1,t_2,t_3]\otimes \Lambda (v_1,v_2 ,v_3),D)\to 
(\Q [t_1,t_2]\otimes \Lambda (v_1,v_2 ,v_3),\overline{D})$$
where
$(\Q [t_1,t_2]\otimes\Lambda (v_1,v_2 ,v_3),\overline{D})=M(ET^2\times_{T^2}X)$ and 
$Dv_1=f_1$,  $Dv_2=f_2$, $Dv_3=f_3$
with $f_1,f_2,f_3$ a regular sequence in $\Q [t_1,t_2,t_3]$ (see Corollary \ref{JL}).
Indeed, then 
$fd(ET^3\times_{T^3}X)=fd(\Q [t_1,t_2,t_3]\otimes \Lambda (v_1,v_2 ,v_3),D)=10$
and $\omega^5\neq 0$ for $\omega =[\lambda_1 t_1+\lambda_2 t_2+\lambda_3 t_3]$
for some $\lambda_i\in \Q$.
Especially,  Proposition \ref{SS}  does not always  deduce
 $r_0(X)=1$
when $X$ is pre-c-symplectic (cf. Theorem \ref{A}).

Recall the {\it Hasse diagram ${\mathcal H}(X)$ of  rational toral ranks} 
for  a simply connected space $X$ \cite{Y}, which 
is the Hasse diagram of a poset induced by 
 ordering 
of the Borel fibrations of rationally almost free toral actions  on $X$.
When 
there exists  a free $t$-toral action on a finite complex $X'$ of same rational homotopy 
with $X$ (Proposition \ref{H}),  
we can describe    a  point $P=[ET^t\times_{T^t}X']$
rationally  presented by the Borel space 
$Y=ET^t\times_{T^t}X'$ 
in the lattice points of the quadrant I.
The coordinate is 
$$P:=(s,t)\ ; \ \ 0\leq s,t,\ \ s+t\leq r_0(X)$$  when 
$ r_0(ET^t\times_{T^t}X')=r_0(X)-s-t$.
In particular,
the root  $(0,0)$
is presented by $X$ itself.
There is an  order $P_i<P_j$  
 given by  the existence of a rational fibration
$$Y_1\to Y_2\to BT^{t_2-t_1}$$
for $P_i=[Y_1]=(s_1,t_1)$ and $P_j=[Y_2]=(s_2,t_2)$
with $s_1\leq s_2$ and $t_1<t_2$.
It is also realized by a $T^{t_2-t_1}$-Borel fibration (Proposition \ref{H}).  
Then $\{ P_i,<\}$ makes a poset and we denote its Hasse diagram
as ${\mathcal H}(X)$.
It may be  useful  to organize knowledge  about almost free toral actions 
(often looks  like the framework of a broken Japanese fan).
  Now, from Proposition   \ref{SS},  we  immediately  obtain
 a necessary condition for  $X$ to be pre-c-symplectic as

\begin{thm} \label{AA}
If $X$ is pre-c-symplectic,  
then  there exists  the  point $P=(r_0(X)-1,1)$ in
${\mathcal H}(X)$.
\end{thm}
It schematically 
gives   a 
necessary condition   
for the existence of a c-symplectic space $Y=ES^1\times_{S^1}X'$
with $X_{\Q}'\simeq X_{\Q}$,  
 in all classes (associated with rational toral ranks) of  orbit spaces of 
rational almost free toral actions on
  $X$.
When $X$ is pre-c-symplectic,
the points  $(r_0(X)-i,i)$
of ${\mathcal H}(X)$, i.e., the leaves of the Hasse  diagram,
may  be presented by c-symplectic models.
For example,
the point $(0,3)$ is surely presented by them
when   $X\simeq_{\Q}S^3\times S^3\times S^7$
as we see in above.
Also see Examples \ref{4} and  \ref{loex}.
When a pre-c-symplectic space $X$ is a product of $n$ odd-spheres,
we can easily check that
there are at least the points
$(2,1), (2,2), \cdots ,(2,n-2)$ in ${\mathcal H}(X)$.
When a  c-symplectic space  is a homogeneous space as in \cite{LO1},
it presents the  point $(0,r_0(X))$
of ${\mathcal H}(X)$ for some pure space $X$ with $\pi_2(X)\otimes \Q=0$
 (see Remark \ref{last}).
On the other hand,
 any  c-symplectic space $Y$ 
 presents   $(r_0(X)-1,1)$
of ${\mathcal H}(X)$ for some pre-c-symplectic  space $X$ with $\dim \pi_2(X)\otimes \Q=
\dim \pi_2(Y)\otimes \Q-1$.

\begin{rem}
The converse of Theorem \ref{AA} is not true.
For example,
put $X=S^3\times  S^3\times S^9\times S^{11}\times S^{13}\times S^{15}\times S^{19}$,
which is not pre-c-symplectic from Theorem \ref{A}
since $k_3+k_4=9+11>19=k_7$ ($n=7$). 
But  there is a point $P=(r_0(X)-1,1)=(6,1)$ in
${\mathcal H}(X)$ presented by 
a model $(\Q [t]\otimes  \Lambda (v_1,..,v_7),D)$
with the differential $Dv_1=\cdots =Dv_4=0$,
$Dv_5=v_2v_3t$, $Dv_6=v_1v_4t$, 
$Dv_7=v_1v_6t+v_2v_5t^2+t^{10}$ in $(4)$
for $H^*(X;\Q)=\Lambda (v_1,..,v_7)$
with $|v_1|=|v_2|=3$,
$|v_3|=9$, $|v_4|=11$, $|v_5|=13$, $|v_6|=15$ and $|v_7|=19$.
We can directly check $r_0(\Q [t]\otimes  \Lambda (v_1,..,v_7),D)=0$
from Proposition \ref{H}.
\end{rem}

This paper is  purely a  Sullivan model  approach to the opening question
restricted on c-symplectic structures 
in the simply connected case. 
Then we see that  the  ratio of degrees in elliptic model structure 
(homotopy rank type \cite{NY})
play an important role to be  pre-c-symplectic.
It consists of three sections.
In \S 2, we give the proof of Theorem \ref{A} 
and  see some related topics.
In particular,  we see in Theorem \ref{odd}
that a space is  pre-c-symplectic
imposes a restrict on the degrees
 when its rational homotopy group is finite oddly generated.
In \S 3, we 
prove Proposition \ref{SS} under a Halperin's criterion 
(Proposition \ref{H})  and 
see some examples of 
 ${\mathcal H}(X)$
when $X$ is pre-c-symplectic in the cases of $r_0(X)\leq 5$.
\\




\noindent
{\bf Acknowledgement}. 
The authors would like to express their gratitude to 
the referee
for his many  valuable comments
to improve the paper.
In particular, he suggested
that they  should rewrite the
introduction to emphasize the toral actions.

\section{Proof   and examples}



In the following 
Lemmas \ref{L1} and \ref{sym}, we assume that $M(X)=(\Lambda (v_1,v_2,\cdots ,v_n),d)$
where
$|v_i|=k_i$ are odd for all $i$ and  $1<k_1\leq \cdots \leq k_n$ for an odd integer $n$.
The symbol $(f_1,..,f_k)$ means the ideal 
of $\Q[t]\otimes \Lambda (v_1,v_2,\cdots ,v_n)$
generated by  elements $f_1,..,f_k$ and `$f\sim g$' means
the $D$-cocycles 
$f$ and $g$ are cohomologuous 
in $(\Q[t]\otimes \Lambda (v_1,v_2,\cdots ,v_n),D)$
of $(2)$; i.e., $[f]=[g]$ in $H^*(Y;\Q )$.


\begin{lem}\label{L1}
If $(\Q[t]\otimes \Lambda (v_1,v_2,\cdots ,v_n),D)$
is c-symplectic, then 
we can put $D$ up to dga-isomorphisms so that\\
${\rm (i)}\ \ \ Dv_i\in (v_1,..,v_{i-1})$
for all $i<n$,\\
${\rm  (ii)}\ \ \ Dv_n=f-\lambda t^{(k_n+1)/2}$
for some $f\in (v_1,v_2,\cdots ,v_{n-1})$ and $\lambda \neq 0\in \Q$,\\ 
${\rm  (iii)}\ \ \ v_1v_2\cdots v_{n-1}\cdot t^{(k_n-1)/2}\sim \lambda t^{(fd(X)-1)/2}$
for some  $\lambda \neq 0\in \Q$.

\end{lem}
\noindent
{\it Proof.}
(i) Suppose that  there is an element $v_i$ with $i<n$
such that $Dv_i=g-\lambda t^{(k_i+1)/2}$ 
for some $g\in (v_1,.., v_{i-1})$ and $\lambda \neq 0\in \Q$.
Then  $\dim H^*(\Q[t]\otimes \Lambda (v_1,v_2,\cdots ,v_i),D)<\infty$
and $fd(\Q[t]\otimes \Lambda (v_1,v_2,\cdots ,v_i),D)=
k_1+\cdots +k_i-1$ \cite{FHT}.
Therefore   we deduce $t^{a/2+1}\sim 0$; i.e., $[t^{a/2+1}]=0$
for some 
 $a<fd (X)-1=k_1+\cdots + k_n-1$.
 It contradicts the definition of a c-symplectic space.
 
 (ii) It is required 
 from (i) and $\dim H^*(\Q[t]\otimes \Lambda (v_1,v_2,\cdots ,v_n),D)<\infty$.
 
(iii)  The element  $v_1v_2\cdots v_{n-1}$ is a $D$-cocycle
from  $Dv_1=Dv_2=0$ and 
 (i).
It is not $D$-exact
from  
(ii).
Then we have 
$[v_1v_2\cdots v_{n-1}]\cdot [t^a]=\lambda [t^{(fd(X)-1)/2}]$
in $H^* (\Q[t]\otimes \Lambda (v_1,\cdots ,v_n),D)$ 
for $a=(fd(X)-1-k_1-\cdots -k_{n-1})/2=(k_n-1)/2$ from the Poincar\'{e} duality property.
\hfill\qed\\

\begin{lem}\label{sym}
Suppose that    $(\Q[t]\otimes \Lambda (v_1,v_2,\cdots ,v_n),D)$
satisfies $Dv_n=f-t^{(|v_n|+1)/2}$
for some $f=g_1t^{a_1}+\cdots +g_kt^{a_k}$ 
with monomials  $g_i\in  \Lambda (v_1,..,v_{n-1})$
and $a_i\geq 0$.
If it 
is c-symplectic,  then  
$g_{i_1}\cdots g_{i_m}\neq 0\in (v_1v_2\cdots v_{n-1})$
for 
some $g_{i_1},..,g_{i_m}$ ($m\leq k$).
\end{lem}

\noindent
{\it Proof.} 
From the assumption, for  $M:=(|v_n|+1)/2$, we have 
$$g_1t^{a_1}+\cdots +g_kt^{a_k}\sim t^M.$$
Suppose $g_{i_1}\cdots g_{i_m}\neq 0$.
By the multiplication  of $t^{M-a_{i_1}}$ on the both sides,
we have $$g_{i_1}g_{i_2}t^{a_{i_2}}+\cdots =g_{i_1}(g_1t^{a_1}+\cdots +g_kt^{a_k})+\cdots\sim g_{i_1}t^{M}+\cdots \  \underset{}{\sim} t^{2M-a_{i_1}}.$$
Again by the multiplication  of $t^{M-a_{i_2}}$ on the both sides,
we have 
$$g_{i_1}g_{i_2}g_{i_3}t^{a_{i_3}}+\cdots  \ \underset{}{\sim} t^{3M-a_{i_1}-a_{i_2}}.$$
Iterate the multiplication  of $t^{M-a_{i_j}}$ to  $j=m-1$.
Then  we have 
$$g_{i_1}g_{i_2}\cdots g_{i_m}t^{a_{i_m}}+\cdots \  \underset{}{\sim} t^{mM-a_{i_1}-\cdots -a_{i_{m-1}}}.$$
Finally  we have
$$g_{i_1}g_{i_2}\cdots g_{i_m}t^{M-1}+\cdots \sim t^{(m+1)M-a_{i_1}-\cdots -a_{i_{m}}-1}=t^{(|g_{i_1}|+\cdots +|g_{i_m}|+|v_n|-1)/2}.$$
If $g_{i_1}\cdots g_{i_m}=\lambda v_1v_2\cdots v_{n-1}$ for some $\lambda\neq 0\in \Q$,
then 
$$ ({\lambda}+\cdots ) v_1v_2\cdots v_{n-1}t^{M-1}\sim  t^{(k_1+k_2+\cdots +k_n-1)/2}=t^{(fd(X)-1)/2}$$
and it  makes a non-zero 
class of $H^{fd(X)-1}(\Q[t]\otimes \Lambda (v_1,v_2,\cdots ,v_n),D)$
when  $\lambda +\cdots \neq 0$.
If there are  no such elements $g_{i_1},g_{i_2},\cdots ,g_{i_m}$,
 then 
$(\Q[t]\otimes \Lambda (v_1,v_2,\cdots ,v_n),D)$
is not c-symplectic  from Lemma \ref{L1}(iii).
\hfill\qed\\
 


\noindent
{\it Proof of Theorem \ref{A}.}
The ``{\it  if}'' part:
We can define the  model 
$(\Q[t]\otimes \Lambda (v_1,v_2,\cdots ,v_n),D)$
of $(2)$ by putting $Dv_1=\cdots =Dv_{n-1}=0$ and 
 $$Dv_n=v_{1}v_{n-1}t^{a_1}+v_{2}v_{n-2}t^{a_2}+\cdots 
+v_{(n-1)/2}v_{(n+1)/2}t^{a_{n-1}}-t^{a_n}$$
for  suitable  $a_i$.
Then $v_{1}v_{n-1}t^{a_1}+v_{2}v_{n-2}t^{a_2}+\cdots 
+v_{(n-1)/2}v_{(n+1)/2}t^{a_{n-1}}\sim t^{a_n}$
deduces, by iterated multiplications  of $t$,
$$v_1\cdots v_{n-1}t^{(k_n-1)/2}\sim t^{(\dim X-1)/2},$$
where    the left  side   is not $D$-exact.
Thus $(\Q[t]\otimes \Lambda (v_1,v_2,\cdots ,v_n),D)$
is c-symplectic

The ``{\it only if}'' part:
From Lemma \ref{L1}(ii),
we can put 
$$Dv_n=\sum_{i=1}^rg_it^{n_i}-t^{(k_n-1)/2}$$ with $g_1,..,g_r$ some  
 monomials  in  $\Lambda (v_1,.., v_{n-1})$
and  $n_i=(|v_n|-|g_i|+1)/2$.
From Lemma \ref{sym},  there is  the set 
$$S:=\{ \  v_{i_1},v_{j_1},  \ \cdots , \  v_{i_{(n-1)/2}},v_{j_{(n-1)/2}} \  \}$$  
such that   
 $S=\{ v_1\cdots ,v_{n-1}\}$ 
and that there are indexes  $l_k$ for  $k=1,..,(n-1)/2$
such that $g_{l_k}$ contains the term  $v_{i_k}v_{j_k}$; i.e.,
$g_{l_k}\in (v_{i_k}v_{j_k})$.
Then
$$|v_{i_k}|+|v_{j_k}|=|v_{i_k}v_{j_k}|\leq |g_{l_k}|<|v_n|$$
for $k=1,..,(n-1)/2$.
From Proposition  \ref{SAT}  below, we have
 $|v_1|+|v_{n-1}|<|v_n|$,
 $|v_2|+|v_{n-2}|<|v_n|$, $\cdots$ and
 $|v_{(n-1)/2}|+|v_{(n+1)/2}|<|v_n|$.\hfill\qed\\



\begin{lem}\label{S}
Let $S=\{a_1,a_2,\dots,a_{2n}\}$ be a set of real numbers with
$a_1 \le a_2 \le \cdots \le a_{2n}$. 
For any partition 
$$
\mathcal{T}=\{\{a_{i_1}, a_{j_1}\},\{a_{i_2}, a_{j_2}\},\dots,\{a_{i_n}, a_{j_n}\}\}
$$
 of $S$ into 2-subsets, where  $i_k,j_k \in \{1,2,\dots,2n\}$ and $i_k \ne j_k$ for $k=1,2,\dots,n$,
 there exsits an element  $\{a_{i_k},a_{j_k}\}$ of $\mathcal{T}$ such that
 $$
 \left\{
 \begin{array}{ccc}
 a_1+a_{2n} &\le& a_{i_k}+a_{j_k} \\
 a_2+a_{2n-1} &\le& a_{i_k}+a_{j_k} \\
 \hdots && \\
 a_n+a_{n+1} &\le& a_{i_k}+a_{j_k}.
 \end{array}
 \right.
 $$
\end{lem}

\noindent
{\it Proof.}
We show the result   by induction on the positive integer $n$. 
For $n=1$, the statement is true since $a_1+a_2 \le a_1+a_2$. Assume the statement is true for $n-1$.
We must prove the assertion is also true for $n$.
Let
$$
\mathcal{T}=\{\{a_{i_1}, a_{j_1}\},\{a_{i_2}, a_{j_2}\},\dots,\{a_{i_n}, a_{j_n}\}\}
$$
be any partition of $S$ into 2-subsets and
let $\{a_i, a_{2n}\} (1\le i \le 2n-1)$ be an element of $\mathcal{T}$ containing $a_{2n}$.

Case of  $a_n \le a_i$.
Then  we have
$$
 \left\{
 \begin{array}{ccc}
 a_1+a_{2n} &\le&a_n+a_{2n} \le  a_i +a_{2n}\\
a_2+a_{2n-1} &\le&a_n+a_{2n} \le  a_i +a_{2n}\\
 \hdots && \\
 a_n+a_{n+1}&\le&a_n+a_{2n} \le  a_i +a_{2n},
 \end{array}
 \right.
 $$
hence we may take $\{a_{i_k},a_{j_k}\}$ as $\{a_i ,a_{2n}\}$.

Case of  $ a_i \le a_{n-1}$.
 Then we have
$$
\hspace{2cm}
\left\{
 \begin{array}{ccc}
 a_1+a_{2n} &\le& a_i +a_{2n}\\
a_2+a_{2n-1} &\le&a_i +a_{2n}\\
 \hdots && \\
 a_i+a_{2n+1-i}&\le&  a_i +a_{2n}. 
 \end{array}
 \right.
\hspace{3cm} (*)
$$
We consider $\mathcal{T'}=\mathcal {T} \backslash \{a_i ,a_{2n}\}$. Since 
$\sharp \mathcal{T'}=n-1$
($\sharp$ denotes the cardinality of a set),
we can apply the induction hypothesis to $\mathcal{T'}$. 
Since $a_1 \le a_2 \le \cdots \le a_{i-1} \le a_{i+1} \le \cdots \le a_{2n-1}$, 
there exsits an element $\{a_{i_k},a_{j_k}\}$ of $\mathcal{T'}$ such that 
$$
\left\{
 \begin{array}{ccc}
 a_1+a_{2n} &\le& a_{i_k}+a_{j_{k}}\\
a_2+a_{2n-1} &\le& a_{i_k}+a_{j_{k}} \\
 \hdots&& \\
 a_{i-1}+a_{2n-i+1} &\le& a_{i_k}+a_{j_{k}} \\
 a_{i+1}+a_{2n-i} &\le& a_{i_k}+a_{j_{k}} \\
  \hdots&& \\
 a_n+a_{n+1}&\le& a_{i_k}+a_{j_{k}}.
 \end{array}
 \right.
 \hspace{3cm} (**)
$$
From $(*)$ and $(**)$, we conclude that
$$
\left\{
 \begin{array}{ccc}
 a_1+a_{2n} &\le& a_{i}+a_{2n}\\
a_2+a_{2n-1} &\le& a_{i}+a_{2n} \\
 \hdots && \\
 a_{i-1}+a_{2n-i+1} &\le& a_{i}+a_{2n} \\
 a_{i+1}+a_{2n-i} &\le& a_{i_k}+a_{j_{k}} \\
  \hdots && \\
 a_n+a_{n+1}&\le& a_{i_k}+a_{j_{k}}.
 \end{array}
 \right.
$$
If we put $Max\{a_i,+a_{2n}, a_{i_k}+a_{j_{k}}\}=a_s+a_t$, 
then $\{a_s,a_t\}$ satisfies the desired inequality.
\hfill\qed\\

From this lemma,  we have   immediately

\begin{prop}{\rm (cf.\cite[Proposition 1.1]{O})}\label{SAT}
Let $S=\{a_1,a_2,\dots,a_{2n}\}$ be a set of positive integers with
$a_1 \le a_2 \le \cdots \le a_{2n}$. 
Assume that there exsits a positive integer $N$ such that
$$
 \left\{
 \begin{array}{ccc}
a_{i_1}+a_{j_1} & \le & N\\
 a_{i_2}+a_{j_2} & \le & N \\
 \hdots && \\
 a_{i_n}+a_{j_n} & \le & N
 \end{array}
 \right.
 $$
for a partition 
$$
\mathcal{T}=\{\{a_{i_1}, a_{j_1}\},\{a_{i_2}, a_{j_2}\},\dots,\{a_{i_n}, a_{j_n}\}\}
$$
 of $S$ into 2-subsets, where  $i_k,j_k\in\{1,2,\dots,2n\}$ and $i_k \ne j_k$ for $k=1,2,\dots,n$.
 Then we have the following inequality:
 $$
 \left\{
 \begin{array}{ccc}
 a_1+a_{2n} &\le& N \\
 a_2+a_{2n-1} &\le& N\\
 \hdots & & \\
 a_n+a_{n+1} &\le & N.
 \end{array}
 \right.
 $$
\end{prop}

 In \cite{O},
we can see
various     versions  of Proposition \ref{SAT}.

From the proof of Lemma \ref{sym}, we have 
\begin{prop}\label{suff}
Suppose that  
 $M(X)= (\Lambda (v_1,v_2,\cdots ,v_n),d)$ with all $|v_i|$ odd  
and that $(\Q[t]\otimes \Lambda (v_1,v_2,\cdots ,v_n),D)$
satisfies $Dv_n=f-t^{(|v_n|+1)/2}$
for some $f=g_1t^{a_1}+\cdots +g_kt^{a_k}$ 
with monomials  $g_j=\lambda_j v_{j_1}\cdots v_{j_{m_j}}\in  \Lambda (v_1,..,v_{n-1})$,
$\lambda_j\neq 0\in \Q$ and $a_j\geq 0$.
If $\prod_{j=1}^k v_{j_1}\cdots v_{j_{m_j}} \neq 0\in (v_1v_2\cdots v_{n-1} )$, then 
it 
is c-symplectic.
\end{prop}

From  the proof of  the ``{\it only if}'' part of Theorem \ref{A}, we have

\begin{thm}\label{odd}
Suppose that
 $M(X)= (\Lambda (v_1,v_2,\cdots ,v_n),d)$ with all $|v_i|$ odd
and $1<|v_1|\leq |v_2|\leq \cdots \leq |v_n|$.
If $X$ is pre-c-symplectic, then $n$ is odd and 
 $|v_1|+|v_{n-1}|\leq |v_n|+1$,
 $|v_2|+|v_{n-2}|\leq |v_n|+1$, $\cdots$,
 $|v_{(n-1)/2}|+|v_{(n+1)/2}|\leq |v_n|+1$.
\end{thm}

\begin{que}
What is the necessary and sufficient condition for a model 
 $(\Lambda (v_1,v_2,$
$\cdots ,v_n),d)$ with
all $|v_i|$ odd
to be pre-c-symplectic ? 
\end{que}

\noindent
{\it Proof of Corollary  \ref{L}.}
The rational types of  compact connected simple Lie groups 
are given as 
$$
\begin{array}{ll}
A_n & (3,5,\dots,2n+1),\\
B_n & (3,7,\dots,4n-1),\\
C_n & (3,7,\dots,4n-1),\\
D_n & (3,7,\dots,4n-5,2n-1),\\
G_2 & (3,11),\\
F_4 & (3,11,15,23), \\
E_6 & (3,9,11,15,17,23),\\
E_7 & (3,11,15,19,23,27,35),\\
E_8 & (3,15,23,27,35,39,47,59)
\end{array}
$$
(see \cite{M}).
For $A_n$, even if $n$ is odd,
we have $3+(2n-1)=2n+1$, which does not 
 satisfy the condition of Theorem \ref{A}.
It is obvious that $B_n$ ($C_n$) and $E_7$ satisfy the condition of Theorem \ref{A}
as
$$3+4(n-1)-1<4n-1,\  7+4(n-2)-1<4n-1, \cdots ,
(2n-3)+(2n+1)<4n-1
$$
$$\mbox{and }\ \ \ 3+27<35,\ \  11+23<35,\ \ 15+19<35,$$
respectively.
Since the ranks of 
$G_2$, $F_4$, $E_6$ and $E_8$ are even,
they are not pre-c-symplectic.
Finally we check $D_n$.
Put   an odd integer $n=2k+1 (k \ge 1)$.
Assume there is an integer $N$ as in Proposition \ref{SAT} for 
the set $S=\{3,7,\dots,8k-5,4k+1\}$.
Then $N=4n-5=4(2k+1)-5=8k-1$.
 Sorting elements of $S$ into increasing order, we have
$$
a_1=3 \le a_2=7 \le \cdots \le a_k=4k-1 \le a_{k+1}=4k+1 \le a_{k+2}=4k+3 \le$$
$$ \cdots \le a_{2k-1}= 8k-9 \le a_{2k}= 8k-5.$$
Then  $a_k+a_{k+1}=(4k-1)+(4k+1)=8k>N$.
It contradicts Proposition \ref{SAT}.
Therefore, Theorem \ref{A}  does not hold for $D_n$.
\hfill\qed\\

\begin{ex}\label{ex2} 
 Even when  a space $X$ is
a product of odd-spheres,
the c-symplectic spaces whose 
pre-c-symplectic space is $X$
are various.
For example, when $X=S^3\times S^5\times S^9\times S^{15}\times S^{33}$,
 there are at least the following twenty rational homotopy types of c-symplectic models $(2)$ with  
the differential $Dv_1=Dv_2=0$ and 
\begin{align*}
1)& \ \ Dv_5=v_1v_4t^8+v_2v_3t^{10}+t^{17}, \ Dv_3=Dv_4=0\\
2)& \ \ Dv_5=v_1v_4t^8+v_2v_3t^{10}+t^{17}, \ Dv_3=0,\ Dv_4=v_1v_2t^4\\
3)& \ \ Dv_5=v_1v_4t^8+v_2v_3t^{10}+t^{17}, \ Dv_3=0,\ Dv_4=v_1v_3t^2\\
4)& \ \ Dv_5=v_1v_4t^8+v_2v_3t^{10}+t^{17}, \ Dv_3=v_1v_2t,\ Dv_4=0\\
5)& \ \ Dv_5=v_1v_4t^{8}+v_2v_3t^{10}+t^{17}, \ Dv_3=v_1v_2t,\ Dv_4=v_1v_3t
\end{align*}
\begin{align*}
6)& \ \ Dv_5=v_1v_2t^{13}+v_3v_4t^5+t^{17}, \ Dv_3=Dv_4=0\\
7)& \ \ Dv_5=v_1v_2t^{13}+v_3v_4t^5+t^{17}, \ Dv_3=0,\ Dv_4=v_1v_3t^2\\
8)& \ \ Dv_5=v_1v_2t^{13}+v_3v_4t^5+t^{17}, \ Dv_3=0,\ Dv_4=v_2v_3t
\end{align*}
\begin{align*}
9)& \ \ Dv_5=v_1v_3t^{11}+v_2v_4t^7+t^{17}, \ Dv_3=Dv_4=0\\
10)& \ \ Dv_5=v_1v_3t^{11}+v_2v_4t^7+t^{17}, \ Dv_3=0,\ Dv_4=v_1v_2t^4\\
11)& \ \ Dv_5=v_1v_3t^{11}+v_2v_4t^7+t^{17}, \ Dv_3=0,\ Dv_4=v_2v_3t\\
12)& \ \ Dv_5=v_1v_3t^{11}+v_2v_4t^7+t^{17}, \ Dv_3=v_1v_2t,\ Dv_4=0\\
13)& \ \ Dv_5=v_1v_3t^{11}+v_2v_4t^7+t^{17}, \ Dv_3=v_1v_2t,\ Dv_4=v_2v_3t
\end{align*}
\begin{align*}
14)& \ \ Dv_5=v_1v_2v_3v_4t+t^{17}, \ Dv_3=Dv_4=0\\
15)&\ \ Dv_5=v_1v_2v_3v_4t+t^{17}, \ Dv_3=0,\ Dv_4=v_1v_2t^4\\
16)& \ \ Dv_5=v_1v_2v_3v_4t+t^{17}, \ Dv_3=0,\ Dv_4=v_1v_3t^2\\
17)& \ \ Dv_5=v_1v_2v_3v_4t+t^{17}, \ Dv_3=0,\ Dv_4=v_2v_3t\\
18)& \ \ Dv_5=v_1v_2v_3v_4t+t^{17}, \ Dv_3=v_1v_2t,\ Dv_4=0\\
19)& \ \ Dv_5=v_1v_2v_3v_4t+t^{17}, \ Dv_3=v_1v_2t,\ Dv_4=v_1v_3t^2\\
20)& \ \ Dv_5=v_1v_2v_3v_4t+t^{17}, \ Dv_3=v_1v_2t,\ Dv_4=v_2v_3t
\end{align*}
for $|v_1|=3,|v_2|=5,|v_3|=9,|v_4|=15,|v_5|=33$.
Note that only $1),\ 6),\ 9)$ and $14)$ are  two stage models  and { formal}; i.e.,
the minimal model is formally constructed 
from its cohomology \cite{LO1}, \cite{FHT}.
Note that   $1)\sim 20)$ 
make a poset structure
as in \cite{YJP}.
For example, 
we have ``$5)\ <\ 3)\ <\ 1)\ <\ 14)\ <\ 0)$'',
where the model $0)$ is given by  $Dv_1=\cdots =Dv_5=0$ (the model of $X$).
 For a product $S^{k_1}\times S^{k_2}\times S^{k_3}\times S^{k_4}\times S^{k_5}$
of odd spheres with $k_1\leq \cdots \leq k_5$,
the inequations
that $$k_1+k_2<k_3,\ \ k_2+k_3<k_4,\ \  k_1+k_2+k_3+k_4<k_5$$
make the most c-symplectic
models.
Conversely,
when  $$k_1+k_2>k_4,\ \ \ k_2+k_4>k_5$$
the c-symplectic model
is uniquely determined up to dga-isomorphism.
For example,
when $(k_1,..,k_5)=(3,5,5,7,11)$,
$$Dv_1=\cdots =Dv_4=0,\ \ Dv_5=v_1v_4t+v_2v_3t+t^6.$$

\end{ex}

\begin{rem}
Put the set C-Symp$(X):=\{$rational  homotopy types of c-symplectic spaces in  $(1)$ with the fibre $X \}$.
Then  C-Symp$(X)=\phi$ if $X$ is not pre-c-symplectic.
For example, 
$\sharp$C-Symp$(S^{k_1}\times S^{k_2} \times S^{k_3})\leq 1$
when $k_i$ are odd, 
$\sharp$C-Symp$(Sp(5))\geq 4$ (see \S 1)
and $\sharp$C-Symp$(S^3\times S^5\times S^9\times S^{15}\times S^{33})\geq 20$
(see Example \ref{ex2}).
When  $Y$ is  c-symplectic and $X$ is pre-c-symplectic,
$Y\times X$ is  pre-c-symplectic
and there is an inclusion 
C-Symp$(X)\subset $ C-Symp$(Y\times X)$ as sets. 
For example, C-Symp$(S^3)=\{ S^2_{\Q} \}$ (one point) and 
 C-Symp$(S^2\times S^3)$ is
$$\{ (\Q [t]\otimes \Lambda (v_1,v_2,v_3),D_a)\ ; \ 
D_av_1=0, D_av_2=tv_1,D_av_3=v_1^2+at^2,\ a\in \Q^* \}/\simeq $$
$\cong \Q^* /{\Q^*}^2$  for $\Q^*:=\Q-0$,  $|v_1|=2$ and $ |v_2|=|v_3|=3$ as a set 
\cite{MS},
which is infinite. 
Also we can give an equivalence 
relation in the rational homotopy types of simply connected
c-symplectic spaces,
that is, put 
$Y\sim Y'$ for two c-symplectic spaces $Y$ and $Y'$ when there are certain finite 
maps
$$Y\leftarrow X_1\to Y_1\leftarrow X_2\to \cdots \to Y_{n-1}\leftarrow X_n\to Y'$$
which  are  fibre inclusions of $(1)$ 
($Y_i$ are c-symplectic).
It satisfies the laws of reflectance, symmetry and transitivity. 
For example,
the models $1),..,20)$ in Example \ref{ex2} are all equivalent.
 \end{rem}
 
\begin{rem}
Recall the rational LS category ${\rm cat}_0(Y)$ of a simply connected space $Y$ 
\cite[{\bf 27}]{FHT}.
It  is equal to the Toomer's invariant of $Y$ (the biggest $s$ for which there
is a non trivial  class in 
$H^*(Y;\Q)=H^*(\Lambda W)$
represented by a cycle in $\Lambda^{\geq s}W$)  
when $Y$ is a rationally Poincar\'{e} duality space(r.P.d.s.) \cite{FHL}.
  For a simply connected space $X$ with $\dim H^*(X;\Q)<\infty$,
put $$c(X)=\mbox{sup} \{ \frac{2{\rm cat}_0(Y)}{fd(X)-1}\ |\ \mbox{fibrations }X\to Y\to K(\Z ,2)\mbox{
where } Y\mbox{ are  r.P.d.s.}\},$$ 
where  $c(X):=0$ if no such  space $Y$ exists  for $X$.
Then $c(X)$ is a rational number with $0\leq c(X)\leq 1$.
In particular,  
i) $c(X)=0$ if $X$ is c-symplectic,
ii)  $c(X)=1$
if and only if $X$ is pre-c-symplectic
and iii) $c(X)\leq c(X\times Y)$ for any c-symplectic space $Y$. 
For example,
when $X_n=S^7\times S^7\times S^{2n+1}$,
$c(X_n)$ is given as 
\begin{center}
{\begin{tabular}{|c||c |c|c|c|c|c|c|c|c|c|}
\hline
$n$&$1$& $2$ & $3$ &  $4$& $5$& $6$ &$7$&$8$&$9$&$\cdots$\\
\hline
$c(X_n)$ &  $\frac{5}{8}$   &$\frac{5}{9}$ &$\frac{1}{2}$&$\frac{6}{11}$&$\frac{7}{12}$&$\frac{8}{13}$&$1$&$1$&$1$&$\cdots$\\
\hline
\end{tabular}
}
\end{center}
When $X_n=S^3\times S^{2n}$,
$c(X_n)=2/(n+1)$ and 
$\lim_n c(X_n)=0$.
When $X_n=S^3\times S^{2n+1}$,
$c(X_n)=(2n+2)/(2n+3)$.
Though $X_n$ is not pre-c-symplectic
for any $n$,
we have 
$\lim_n c(X_n)=1$.
\end{rem}
 
\begin{ex} \label{cpn}
For any product of odd-spheres
$X=S^{k_1}\times \cdots \times S^{k_n}$
with $n$ odd and $k_1\leq \cdots \leq k_n$,
the product 
$X\times \C P^N$ is pre-c-symplectic
if
$k_1+k_{n-1}\leq 2N$, $k_2+k_{n-2}\leq 2N$,
$\cdots $, $k_{(n-1)/2}+k_{(n+1)/2}\leq 2N$
and $k_n\leq 2N+1$.
Indeed, we can put
 $Dx=Dv_1=\cdots =Dv_{n-1}=0$, $Dv_n=x^{(k_n-1)/2}t$
and $$Dy=x^{N+1}+v_1v_{n-1}t^*+\cdots
+v_{(n-1)/2}v_{(n+1)/2}t^*+t^{N+1} $$
for $M(\C P^N)=(\Lambda (x,y),d)$
with $|x|=2$, $dx=0$ and $dy=x^{N+1}$.
Then $[t^a]\neq 0$ for $a=(k_1+..+k_n-1)/2+N$.
\end{ex}

\begin{rem}
What additional properties  of a c-symplectic space $Y$ (or model $M(Y)$)
can be deduced from the pre-c-symplectic space $X$ in $(1)$ ?
A  c-symplectic space $Y$
of $fd(Y)=2m$  is said that it satisfies  the {hard Lefschetz condition}
with respect to the c-symplectic class $t$ when  the maps
$$\cup  t^k:H^{m-k}(Y;\Q)\to H^{m+k}(Y;\Q)\ \ \ \ 1\leq k\leq m$$
are isomorphisms
\cite{TO}.
For example, a compact K\"{a}hler manifold satisfies  the { hard Lefschetz condition} \cite{TO},  \cite[Theorem 4.35]{FOT}.
As well as when $(\Q [t]\otimes \Lambda V,D)$ of  $(2)$
is c-symplectic, whether or not  it satisfies  
the  hard Lefschetz condition depends on  $D$.
For example, when $H^*(X;\Q)=\Lambda (v_1,v_2,v_3,v_4 ,v_5)$
with $|v_1|=|v_2|=3$, $|v_3|=|v_4|=5$ and $|v_5|=11$,  
put  $Dv_1=\cdots =Dv_4=0$
and 
  $$a)\ \ \ Dv_5=v_1v_2t^3+v_3v_4t+t^6$$
$$b) \ \ \ Dv_5=v_1v_4t^2+v_2v_3t^2+t^6 ,$$
which are both c-symplectic with $m=13$. 
Then $a)$ satisfies 
 the  hard Lefschetz condition  but $b)$ does  not. 
Indeed, \\
Case of  $a)$. When  $k=10$,
$ Ker (\cup  t^{10}:H^{3}(Y;\Q)\to H^{23}(Y;\Q))=0$
since $[v_1t^{10}]=-[v_1(v_1v_2t^3+v_3v_4t)t^4]=-[v_1v_3v_4t^5]\neq 0.$
When  $k=8$,
$ Ker (\cup  t^{8}:H^{5}(Y;\Q)\to H^{21}(Y;\Q))=0$
since $[v_3t^{8}]=-[v_3(v_1v_2t^3+v_3v_4t)t^2]=-[v_1v_2v_3t^5]\neq 0.$
When $k\neq 8,\ 10$,
we can easily check $ Ker (\cup  t^{k})=0$.\\
Case of  $b)$. When  $k=10$, $Ker (\cup  t^{10}:H^{3}(Y;\Q)\to H^{23}(Y;\Q))\neq 0$.
Indeed, 
$[v_1]\in Ker (\cup  t^{10})$
since $$[v_1t^{10}]=-[v_1(v_1v_4t^2+v_2v_3t^2)t^4]=-[v_1v_2v_3t^6]=[v_1v_2v_3(v_1v_4t^2+v_2v_3t^2)]=0.$$


  \end{rem}

\begin{rem}
When a map $g:(Y_1,w_1)\to (Y_2,w_2)$ between simply connected 
c-symplectic spaces 
induces  $H^*(g)(w_2)=w_1$; i.e.,  a {\it c-symplectic map},
  there is a map between fibrations:
$$
\xymatrix{
X_1\ar[r]\ar[d]_f& Y_1\ar[r]\ar[d]^g& K(\Z ,2)\ar@{=}[d]\\
X_2\ar[r]& Y_2\ar[r]& K(\Z ,2),\\
}$$
where  $f:X_1\to X_2$ is the induced map between pre-c-symplectic spaces.
Conversely,  when is  a map $f:X_1\to X_2$ between pre-c-symplectic spaces
extended  to   a c-symplectic map; i.e.,  a {\it pre-c-symplectic map} ?
Refer  \cite{SY} in  the  case of self homotopy equivalences.
  \end{rem}
  
%

\section{Rational toral ranks}

If an $r$-torus $T^r$ acts on a simply connected space $X$
by $\mu :T^r\times X\to X$, there is the Borel fibration
$$
X \to ET^r \times_{T^r} X \to BT^r,
$$
where 
$ ET^r \times_{T^r} X $  is the orbit space of the  action
$g(e,x)=(e\cdot g^{-1},g\cdot x)$  
on the product $ ET^r \times  X $ for $g\in T^r$.
Note that $ET^r \times_{T^r} X$ 
is rational  homotopy equivalent  to the $T^r$-orbit space of $X$
when $\mu$ is an  almost free toral action \cite{FOT}.
The above Borel fibration  is rationally given by the  KS model
$$
(\Q[t_1,\dots,t_r],0)
 \to (\Q[t_1,\dots,t_r] \otimes \Lambda {V},D)
 \to (\Lambda {V},d)\ \ \ \ (4)$$
where    with $|{t_i}|=2$ for $i=1,\dots,r$, $Dt_i=0$ and
$Dv \equiv dv$ modulo the ideal $(t_1,\dots,t_r)$ for $v\in V$.
It is a generalization of $(2)$.
Recall 
Halperin's

\begin{prop}\cite[Proposition 4.2]{H}\label{H}
Suppose that $X$ is a simply connected CW-complex  with 
$\dim H^*(X;\Q)<\infty$.
Put $M(X)=(\Lambda V,d)$.
Then  $r_0(X) \ge r$ if and only if there is a KS model  $(4)$
 satisfying $\dim H^*(\Q[t_1,\dots,t_r] \otimes \Lambda {V},D)<\infty$.
 Moreover,
 if  $r_0(X) \ge r$,
 then $T^r$ acts freely on a finite complex $X'$
 that has  the same rational homotopy type as $X$
 and $M(ET^r\times_{T^r}X')\cong 
(\Q[t_1,\dots,t_r] \otimes \Lambda {V},D)$.
\end{prop}

\noindent
{\it Proof of Proposition   \ref{SS}.} 
Put the formal dimension of $Y$  as  $2n$.
Then there is an element $[\omega]\in H^2(Y;\Q)$ with  $[\omega]^n\neq 0$. 
Suppose  $r_0(Y)>0$.
From Proposition \ref{H},
there is a finite complex $Y'$ with $Y'_{\Q}\simeq Y_{\Q}$ 
and there is a free $S^1$-action on $Y'$.
Thus we have the Borel fibration 
$Y'\overset{i}{\to}  ES^1\times_{S^1}Y'\to BS^1$,
where $[\omega ]$ is a restriction of an element $[u]$ of $H^2(ES^1\times_{S^1}Y';\Q)$;
i.e.,
$i^*([u])=[w]$.
Since  the formal dimension of $ES^1\times_{S^1}Y'$ is $2n-1$,
we have 
$[u]^n=0$.
This is a contradiction.
\hfill\qed\\

Recall the following proposition induced by \cite[Lemma 2.12]{JL}.

\begin{prop}\cite[Lemma 2.1]{Y3}\label{AAAA}  When $X$ is the product of 
$n$ odd-spheres,
the second row of 
${\mathcal H}(X)$ is empty, that is,
there is no point $P=(1,*)$ in  ${\mathcal H}(X)$ for $*=1,2,..,n-1$.
\end{prop}


\begin{cor}\label{JL} For a fibration  
$S^{k_1}\times \cdots \times S^{k_n}
\to X\to \C P^{\infty}\times \cdots \times \C P^{\infty}$
($n-1$-factors)
with $k_1,..,k_n$ odd,
$X$ is pre-c-symplectic
if $\dim H^*(X;\Q)<\infty$.
\end{cor}

\noindent
{\it Proof.} Put $M(S^{k_1}\times \cdots \times S^{k_n})=(\Lambda (v_1,\cdots ,v_n),0)$.
We show that the   model  $M(X)=(\Q [t_1,..,t_{n-1}]\otimes \Lambda (v_1,\cdots ,v_n),D)$ 
is pre-c-symplectic. 
From Proposition \ref{AAAA}(\cite[Lemma 2.12]{JL}),
there is a KS model $(2)$   
$$(\Q [t_n],0)\to (\Q [t_1,..,t_{n}]\otimes \Lambda (v_1,\cdots ,v_n),D')\to (\Q [t_1,..,t_{n-1}]\otimes \Lambda (v_1,\cdots ,v_n),D)$$
 such that the formal dimension 
of $B:=(\Q [t_1,..,t_n]\otimes \Lambda (v_1,\cdots ,v_n),D')$
is $N:=|v_1|+\cdots +|v_n|-n$.
It is formal and the cohomology algebra is 
$$\Q [t_1,\cdots ,t_n]/(D'v_1,\cdots , D'v_n)$$
where $D'v_1,\cdots , D'v_n$ is a regular sequence
in $\Q [t_1,\cdots ,t_n]$.
Then 
$(\lambda_1 t_1+\cdots +\lambda_n t_{n})^{N/2}$
is the fundamental  class of $H^*(B)$ for 
an element $\lambda_1 t_1+\cdots +\lambda_n t_{n}
\in H^2(B)$
with $\lambda_i\in \Q$.
\hfill\qed\\

Thus, when $X$ is a  product of 
$n$ odd-spheres,
 the point $(0,n-1)$ in  ${\mathcal H}(X)$
 is surely  presented by pre-c-symplectic models and 
the point $(0,n)$
 is  by c-symplectic models.
 In the following examples,
 $P_0=(0,0)=[X]$.

\begin{ex}
For a  pre-c-symplectic space $X$ with $r_0(X)=1$, 
the  Hasse diagram ${\mathcal H}(X)$ is (uniquely) given as  
$${\small
\xymatrix{ 
P_1&\\
P_0\ar@{-}[u]&
}}$$
where the point $P_1$ is 
presented by a c-symplectic model.
For example, 
when $X=S^{2n+1}$,
 $P_1=(0,1)=[\C P^n]$.

When $M(X)=(\Lambda (v_1,..,v_{2n+1}),d)$
with $$dv_i=0\ (i<2n+1),\ \ \ dv_{2n+1}=v_1\cdots v_{2j_1}+\cdots +v_{2j_{k-1}+1}\cdots v_{2j_k}\ \ (2j_k=2n),$$
we can put $Dv_i=0$ for $i\neq 2n+1$ and 
$$Dv_{2n+1}=v_1\cdots v_{2j_1}+\cdots +v_{2j_{k-1}+1}\cdots v_{2j_k}+t^{|v_{2n+1}|+1/2}.$$
Then it is formal and c-symplectic from Proposition \ref{suff}.

When $M(X)=(\Lambda (v_1,..,v_{n}),d)$
with $|v_1|=|v_2|=3$, $|v_3|=5$, $\cdots$, $|v_n|=2n-1$ and 
$$dv_1=dv_2=0,\ dv_3=v_1v_2,\  dv_4=v_1v_3,\ \cdots , \  dv_n=v_1v_{n-1}$$
for an odd integer $n>2$,
we can put $Dv_i=dv_i$ for $i\neq n$ and 
$$Dv_n=v_1v_{n-1}+v_2v_{n-2}t-v_3v_{n-3}t+\cdots +(-1)^{a}v_{a}v_{a+1}t+t^{n}$$
for $a=(n-1)/2$.
Then $D\circ D=0$  and it is c-symplectic from Proposition \ref{suff}.
\end{ex}

\begin{ex}
For a  pre-c-symplectic space $X$ with $r_0(X)=2$, 
the  Hasse diagram ${\mathcal H}(X)$ is uniquely given as  
$${\small
\xymatrix{ 
P_2 & \\
P_1 \ar@{-}[u] & P_3\\
P_0\ar@{-}[u]\ar@{-}[ur]&
}}$$
, which has the point $P_3=(1,1)$
from Theorem \ref{AA}.
For example,
it is given when $M(X)=(\Lambda (v_1,v_2,v_3,v_4,v_5),d)$
where $dv_1=dv_2=dv_3=0$, $dv_4=v_1v_2$ and  $dv_5=v_1v_3$
with 
 $|v_1|=|v_2|=3$, $|v_3|=7$, $|v_4|=5$, $|v_5|=9$.
Then 
$P_2=(0,2)=[(\Q [t_1,t_2]\otimes \Lambda (v_1,v_2 ,v_3,v_4,v_5),D)]$ 
where  
 $Dv_1=Dv_2=Dv_3=0$,  
$Dv_4=v_1v_2+t_1^{3}$ 
and $Dv_5=v_1v_3+t_2^{5}$.
 Also $P_3=[(\Q [t]\otimes \Lambda (v_1,v_2 ,v_3,v_4,v_5),D)]$ 
where  
 $Dv_1=Dv_2=Dv_3=0$, $Dv_4=v_1v_2$ and  
$Dv_5=v_1v_3+v_2v_4t+t^{5}$,
which is c-symplectic from Proposition \ref{suff}.
Indeed, $[t^{13}]=[v_1v_2v_3v_4t^4]\neq 0$.
\end{ex}

\begin{ex}\label{3}(see  \cite[Examples 3.5, 3.6]{Y})
Suppose that  $X$ with $r_0(X)=3$ is pre-c-symplectic.
When  $X=S^{k_1}\times S^{k_2}\times S^{k_3}$,
from Theorem \ref{AA} and Proposition \ref{AAAA},
the  Hasse diagram ${\mathcal H}(X)$ is uniquely given as  
$${\small
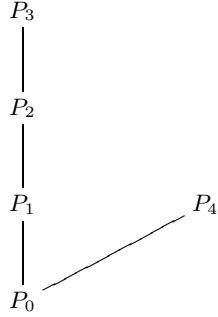
\xymatrix{ 
P_3 &  & \\
P_2 \ar@{-}[u] && \\
P_1 \ar@{-}[u] && P_{4}\\
P_0\ar@{-}[u]\ar@{-}[urr]&&
}}$$
, which has the point $P_4=(2,1)$.
For example,
 when $(k_1,k_2,k_3)=(3,3,7)$, 
 $P_1=[S^2\times S^3\times S^7]$,
$P_2=[S^2\times S^2\times S^7]$
and $P_3=[S^2\times S^2\times  \C P^3]$. 
Here $P_4=(2,1)=[Y]$ is given 
  by 
the  model  $M(Y)=(\Q [t]\otimes \Lambda (v_1,v_2 ,v_3),D)$
with $Dv_1=Dv_2=0$ and  
$Dv_3=v_1v_2t+t^{4}$,
which is c-symplectic.

Next put $M(X)=(\Lambda V,d)=(\Lambda (v_1,v_2,v_3,v_4,v_5),d)$
with $dv_1=dv_2=dv_4=dv_5=0$
and $dv_3=v_1v_2$. 
If 
 $|v_1|=|v_2|=3$, $|v_3|=5$, $|v_4|=9$
and  $|v_5|=13$, 
then
 ${\mathcal H}(X)$ is  given as
$$ {\small 
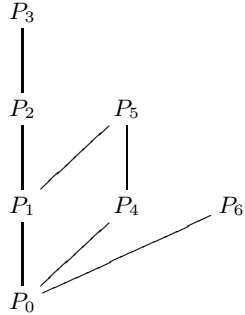
\xymatrix{ 
P_3 &  &   \\
P_2 \ar@{-}[u]& P_5&  \\
P_1 \ar@{-}[u]\ar@{-}[ru] & P_4\ar@{-}[u]& P_6 \\
P_0\ar@{-}[u]\ar@{-}[ur]\ar@{-}[urr]\\
}}
$$
, where 
  $P_3=[ (\Q [t_1,t_2,t_3]\otimes \Lambda V,D)]$ 
with  $Dv_3=v_1v_2+t_2^3$,
$Dv_4=t_1^5$, $Dv_5=t_3^7$, 
$P_4=[ (\Q [t_1]\otimes \Lambda V,D)]$ 
with    $Dv_3=v_1v_2$,
$Dv_4=v_1v_3t_1+t_1^5$, $Dv_5=0$,
 $P_5=[ (\Q [t_1,t_2]\otimes \Lambda V,D)]$ 
with   $Dv_3=v_1v_2$, 
$Dv_4=v_1v_3t_1+t^5_1$, 
$Dv_5=t^7_2$ and 
$P_6=[ (\Q [t]\otimes \Lambda V,D)]$ 
with  $Dv_4=0$,   $Dv_3=v_1v_2$, 
$Dv_5=v_2v_4t+v_1v_3t^3+t^7$.
Here $Dv_1=Dv_2=0$ for all.
This  model presenting  $P_6=(2,1)$
makes $X$ to be pre-c-symplectic
from Proposition \ref{suff}.
Indeed,
$[t^{16}]= [v_1v_2v_3v_4t^6]\neq 0$
for $fd(\Q [t]\otimes \Lambda V,D)=32$.

If $|v_1|=|v_2|=3$, $|v_3|=5$, $|v_4|=9$ and $|v_5|=11$,
it satisfies the necessary condition of Theorem \ref{odd}
that $3+9\leq 11+1$ and $3+5\leq 11+1$.
But we can easily  check that there is no point $P_6=(2,1)$
since 
$Dv_5\in (t,v_1,v_2,v_3)$
in any dga $(\Q [t]\otimes \Lambda V,D)$
from degree reason. 
Indeed,
then $r_0(\Q [t]\otimes \Lambda V,D)>0$ 
since we can put
$D_2(v_4)=t_2^5$ and $D_2(v_i)=D(v_i)$ for $i\neq 4$
as a relative model of $(4)$
$$
(\Q[t_2],0)
 \to (\Q[t_2,t] \otimes \Lambda {V},D_2)
 \to (\Q [t]\otimes \Lambda {V},D)$$
 with $\dim H^*(\Q[t_2,t] \otimes \Lambda {V},D_2)<\infty$.
 Thus ${\mathcal H}(X)$ is  given as
$$ {\small 
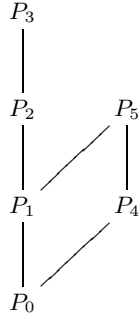
\xymatrix{ 
P_3 &  &   \\
P_2 \ar@{-}[u]& P_5&  \\
P_1 \ar@{-}[u]\ar@{-}[ru] & P_4\ar@{-}[u]&  \\
P_0\ar@{-}[u]\ar@{-}[ur]\\
}}
$$
and $X$ is  not pre-c-symplectic
from Theorem \ref{AA}.
\end{ex}

\begin{ex}\label{4}
Put $M(X)=(\Lambda (v_1,v_2,v_3,v_4,v_5,v_6,v_7),d)$
with $dv_1=dv_2=dv_3=dv_4=dv_7=0$, $dv_5=v_1v_2$, $dv_6=v_1v_3$ and 
$|v_1|=|v_2|=|v_3|=3$,  $|v_4|=|v_5|=|v_6|=5$, $|v_7|=9$. Then
$r_0(X)=4$ and 
 ${\mathcal H}(X)$ is  given as 

$${\small 
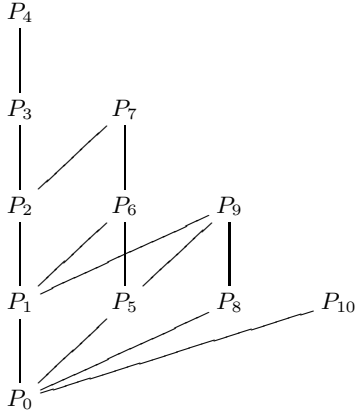
\xymatrix{ 
P_4 & &  & \\
P_3 \ar@{-}[u]& P_7 &  & \\
P_2 \ar@{-}[u]\ar@{-}[ru]& P_6 \ar@{-}[u]& P_{9}& \\
P_1 \ar@{-}[u]\ar@{-}[ru]\ar@{-}[urr] & P_5\ar@{-}[u]\ar@{-}[ru] & P_8\ar@{-}[u]& P_{10} \\
P_0\ar@{-}[u]\ar@{-}[ur]\ar@{-}[urr]\ar@{-}[urrr]\\
}}$$
, where 
the edge
$P_5P_9$ ($P_5<P_9$)
is given by $Dv_i=dv_i$ for $i\neq 4,7$, 
$$Dv_7=v_1v_6t_1+v_2v_5t_2+t_1^5, \ \ Dv_4=t_2^3$$
and 
$P_{10}=(3,1)$ is presented  by $Dv_i=dv_i$ for $i\neq 7$, 
$$Dv_7=v_1v_6t+v_2v_5t+v_3v_4t+t^5,$$
which is c-symplectic from Proposition \ref{suff}.
Also $P_7$ is presented  by a c-symplectic model  with $Dv_i=dv_i$ for $i=1,2,3$, 
$$Dv_7=v_1v_6t_i+t_i^5,\ Dv_5=v_1v_2+t_j^3,\ Dv_4=t_k^3,$$
which gives the sequence of orders
$P_0<P_5<P_6<P_7$ when $(i,j,k)=(1,2,3)$.
\end{ex}

\begin{ex}\label{loex}
When the product of five odd-spheres  $X=S^{k_1}\times S^{k_2}\times  S^{k_3}\times S^{k_4}\times  S^{k_5}$ 
is pre-c-symplectic,
there are (at least) the  following  two   Hasse diagrams $(a)$ and $(b)$ that 
have the point $P_9=(4,1)$.
$${\small
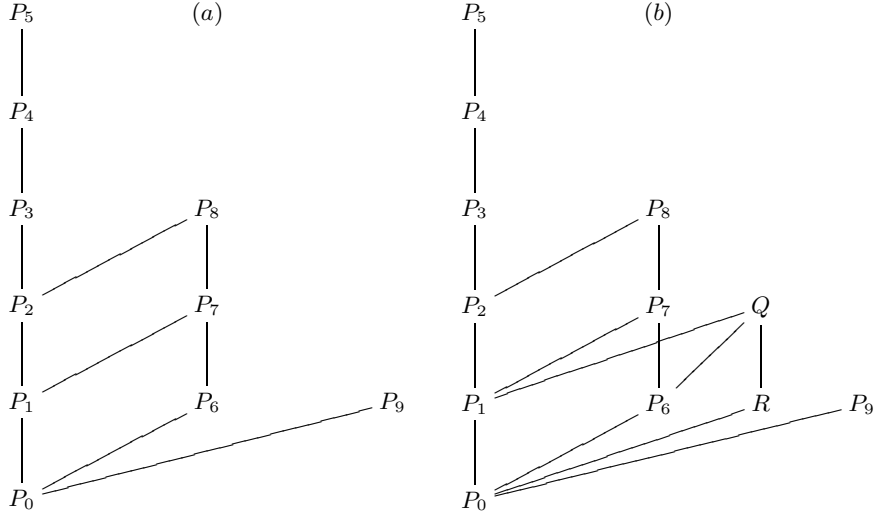
\xymatrix{ 
P_5&&(a)&&\\
P_4 \ar@{-}[u]& &  & &\\
P_3 \ar@{-}[u]&  &P_8&  & \\
P_2 \ar@{-}[u]\ar@{-}[urr]&  & P_{7}\ar@{-}[u]&& \\
P_1 \ar@{-}[u]\ar@{-}[urr] & & P_6\ar@{-}[u]& &P_{9} \\
P_0\ar@{-}[u]\ar@{-}[urr]\ar@{-}[urrrr]\\
}\ \ \ \ \ 
\xymatrix{ 
P_5&&(b)&&\\
P_4 \ar@{-}[u]& &  & &\\
P_3 \ar@{-}[u]&  &P_8&  & \\
P_2 \ar@{-}[u]\ar@{-}[urr]&  & P_{7}\ar@{-}[u]&Q& \\
P_1 \ar@{-}[u]\ar@{-}[urr] \ar@{-}[urrr]& & P_6\ar@{-}[u]\ar@{-}[ur]& R\ar@{-}[u] &P_{9} \\
P_0\ar@{-}[u]\ar@{-}[urrr]\ar@{-}[urr]\ar@{-}[urrrr]\\
}
}$$
For example, $(a)$ is given when $X=S^3\times S^3\times S^3\times S^3\times S^{9}$
and $(b)$ is given when $X=S^3\times S^3\times S^7\times S^{11}\times S^{15}$.
They satisy the condition of Theorem \ref{A}.
The point  $R$ of $(b)$  is presented by the model, for example,  with
 $Dv_1=Dv_2=Dv_5=0$,
$Dv_3=v_1v_2t_1$ and 
$Dv_4=v_1v_3t_1+t_1^{6}$. 
The point  $Q$ of $(b)$  is presented by the model, for example,  with
  $Dv_1=Dv_2=0$,
$Dv_3=v_1v_2t_1$, 
$Dv_4=v_1v_3t_1+t^{6}_1$ and 
$Dv_5=t^8_2$. 
The points  $P_6$  of $(a),(b)$ are  presented by the model, for example,   with $Dv_1=Dv_2=Dv_3=Dv_4=0$ 
and $Dv_5=v_1v_4t^{(k_5-k_1-k_4+1)/2}+t^{(k_5-1)/2}$. 
Finally,  the points  $P_9$  of $(a),(b)$ are  presented by the model, for example,   $Dv_1=Dv_2=Dv_3=Dv_4=0$, 
 $(a):\ Dv_5=v_1v_4t^2+v_2v_3t^2+t^5$
and  $(b):\ Dv_5=v_1v_4t+v_2v_3t^3+t^8$,
which are c-symplectic models.
In these examples of $X$, three points $P_5$, $P_8$ and $P_9$ are presented by
c-symplectic models,  in $(a)$ and $(b)$.
In particular, 
for $M(S^3\times S^3\times S^3\times S^3\times S^{9})=(\Lambda {V},0)$
giving  $(a)$,
the c-symplectic model $ (\Q [t_1,t_2,t_3]\otimes \Lambda {V},D)$
with $(*)$
$$Dv_1=Dv_2=0, \ Dv_3=t_i^2, \ Dv_4=t_j^2, \ Dv_5=v_1v_2t_k^2+t_k^5,$$
where $\{ i,j,k\} =\{ 1,2,3\}$, 
presents $P_8$ and its process of  fibrations
gives the sequence of orders
$P_0<P_1<P_2<P_8$,
$P_0<P_1<P_7<P_8$
or   $P_0<P_6<P_7<P_8$.
On the other hands,    the 
 c-symplectic model $ (\Q [t_1,t_2,t_3]\otimes \Lambda {V},D)$ 
of  Lupton-Oprea\cite[Example 2.12]{LO1} with $(**)$
$$Dv_1=t_i^2, Dv_2=t_it_j, Dv_3=t_j^2, Dv_4=t_jt_k, Dv_5=
t_k^5+(v_1t_j-t_iv_2)(v_3t_k-t_jv_4)$$
presents $P_8$  but  
can  not  give    $P_0<P_6<P_7<P_8$,
especially since
$v_1t_1^2v_4=\overline{D}(-v_1v_3v_4)$
in $(\Q [t_1]\otimes \Lambda V,\overline{D})$
when  $j=1$.
Notice  that the model of $(*)$ is formal  but $(**)$ is not. 

\end{ex}

\begin{rem}\label{last}
Simply connected c-symplectic spaces $Y$
 are schematically  classified by  the following diagrams  ${\mathcal P}(Y)$
 with respect to 
rational toral ranks.
When $\dim \pi_2(Y)\otimes \Q=n$
with $M(Y)=(\Lambda U,d_U )$,
there is the  relative model 
$$(\Q [t_1,..,t_n],0)\to (\Lambda U,d_U )\to (\Lambda V,d)\ \ ; \ V^2=0$$
with $|t_i|=2$ and $U=V\oplus \Q (t_1,..,t_n)$.
Then $Y$ presents a point (leaf)  in 
${\mathcal H}( \Lambda V,d)$
with certain sequences
$[(\Lambda V,d)]<\cdots <[Y]$ of  orders
which are  given by compositions of  fibrations.
Glue all such paths $[(\Lambda V,d)]-\cdots -[Y]$
from $[(\Lambda V,d)]$ to $[Y]$
 in ${\mathcal H}( \Lambda V,d)$
and denote  it as ${\mathcal P}(Y)$.
For example, in the case of $n=3$,
we can find the following  four types of ${\mathcal P}(Y)$ in this paper:
{\small
$${\small
\xymatrix{ 
\bullet   & \\
\bullet \ar@{-}[u] & \\
\bullet \ar@{-}[u] &\\
\bullet\ar@{-}[u]&
}}
{\small
\xymatrix{ 
   & \bullet\\
\bullet \ar@{-}[ur] &\bullet  \ar@{-}[u]\\
\bullet \ar@{-}[u] \ar@{-}[ur]&\bullet \ar@{-}[u]\\
\bullet\ar@{-}[u]\ar@{-}[ur]&
}}\ \ \ \ \ \ 
{\small
\xymatrix{ 
   & &\bullet\\
\bullet \ar@{-}[urr] &&\bullet  \ar@{-}[u]\\
\bullet \ar@{-}[u] \ar@{-}[urr]&&\bullet \ar@{-}[u]\\
\bullet\ar@{-}[u]\ar@{-}[urr]&&
}}\ \ \ \ \ \ 
{\small
\xymatrix{ 
   & &\bullet\\
 \bullet \ar@{-}[urr]&&\bullet  \ar@{-}[u]\\
\bullet \ar@{-}[u]  \ar@{-}[urr]&&\\
\bullet\ar@{-}[u]&&
}}
$$}
, which are in  Example \ref{3}, Example \ref{4}, Example \ref{loex}$(a)(*)$
and  Example \ref{loex}$(a)(**)$, 
 respectively.
 If a c-symplectic space is a homogeneous
 space,
 it is the first type from $r_0(X)\leq -\chi_{\pi}(X):=\dim \pi_{odd}(X)\otimes \Q-\dim \pi_{even}(X)\otimes \Q$
 for an elliptic space  $X$ (\cite{AH})  and 
 \cite[Corollary 2.3]{LO1}.
\end{rem}

\end{document}